\documentclass{article}

\usepackage[english]{babel}

\usepackage[letterpaper,top=2cm,bottom=2cm,left=3cm,right=3cm,marginparwidth=1.75cm]{geometry}

\usepackage{amsmath}
\usepackage{amssymb}
\usepackage{graphicx}
\usepackage{subcaption}
\usepackage{natbib}

\usepackage{times}
\usepackage[utf8]{inputenc}
\usepackage[T1]{fontenc}
\usepackage{url}
\usepackage{booktabs}
\usepackage{multirow}
\usepackage{amsfonts}
\usepackage{nicefrac}
\usepackage{microtype}
\usepackage{faktor}
\usepackage{hhline}

\usepackage{thm-restate}

\usepackage{tikz}
\usepackage{mathtools, stmaryrd}
\usepackage{tocloft}
\usepackage{enumitem}
\usepackage[hidelinks]{hyperref}
\usepackage{xcolor}
\hypersetup{
    colorlinks,
    linkcolor={red!50!black},
    citecolor={blue!50!black},
    urlcolor={blue!80!black}
}

\usepackage[algo2e, algoruled, boxed, vlined]{algorithm2e}
\usepackage{algorithm}
\usepackage{algpseudocode}

\SetKwInput{KwInput}{Input}
\SetKwInput{KwOutput}{Output}

\usepackage{cleveref}

\newtheorem{Th}{Theorem}[section]

\newtheorem{Rem}[Th]{Remark}
\newtheorem{?}[Th]{Problem}

\def\Sp{\operatorname{Sp}}

\usepackage{pifont}

\usepackage{float}
\usepackage{wrapfig}

\newcommand{\Span}{\mathrm{span}}
\newcommand{\argmin}{\arg\!\min}

\newcommand{\xs}{x_\star}

\title{Quadratic minimization: from conjugate gradient to an adaptive Heavy-ball method with Polyak step-sizes}
\author{
Baptiste Goujaud\footnote{CMAP, École Polytechnique, Institut Polytechnique de Paris.
baptiste.goujaud@polytechnique.edu.},
Adrien Taylor\footnote{INRIA, École Normale Supérieure, CNRS, PSL Research University, Paris.
adrien.taylor@inria.fr.},
Aymeric Dieuleveut\footnote{CMAP, École Polytechnique, Institut Polytechnique de Paris.
aymeric.dieuleveut@polytechnique.edu.}
}

\renewcommand{\leq}{\leqslant}
\renewcommand{\geq}{\geqslant}

\renewcommand{\preceq}{\preccurlyeq}
\begin{document}
\addtocontents{toc}{\protect\setcounter{tocdepth}{0}}
\maketitle

\begin{abstract}
    In this work, we propose an adaptive variation on the classical Heavy-ball method for convex quadratic minimization.
    The adaptivity crucially relies on so-called ``Polyak step-sizes'',
    which consists in using the knowledge of the optimal value of the optimization problem at hand
    instead of problem parameters such as a few eigenvalues of the Hessian of the problem.
    This method happens to  also be equivalent to a variation of the classical conjugate gradient method,
    and thereby inherits many of its attractive features, including its finite-time convergence,
    instance optimality, and its worst-case convergence rates.

    The classical gradient method with Polyak step-sizes is known to behave very well in situations in which it can be used,
    and the question of whether incorporating momentum in this method is possible and can improve the method itself appeared to be open.
    We provide a definitive answer to this question for minimizing convex quadratic functions,
    a arguably necessary first step for developing such methods in more general setups.
\end{abstract}

\section{Introduction}\label{sec:introduction}
    Consider the convex quadratic minimization problem in the form
    \begin{equation}
        \min_{x\in\mathbb{R}^d} \left\{f(x)\triangleq \frac{1}{2} \langle x,\, H x \rangle + \langle h,\, x \rangle \triangleq \frac{1}{2} \langle x-\xs,\, H (x-\xs)\rangle + f_\star\right\} \label{eq:main-opt}
    \end{equation}
    where $H\succcurlyeq 0$ is a symmetric positive semi-definite matrix, and we denote $f_\star$ the minimum value of $f$ (a few instances of such problems are presented in~\Cref{sec:experiments}).
    In the context of large-scale optimization (i.e.\ $d\gg 1$), we are often interested in using first-order iterative methods for solving~\cref{eq:main-opt}.
    There are many known and celebrated iterative methods for solving such quadratic problems, including conjugate gradient, Heavy-ball methods (a.k.a., Polyak momentum), Chebyshev methods, and gradient descent.
    Each of those methods having different specifications, the choice of the method largely depends on the application at hand.
    In particular, a typical drawback of momentum-based methods is that they generally require the knowledge of some problem parameters (such as extreme values of the spectrum of $H$).
    This problem is typically not as critical for simpler gradient descent schemes with no momentum, although it generally still requires some knowledge on problem parameters if we want to avoid using linesearch-based strategies.
    This limitation motivates the search for adaptive strategies, fixing step-size using past observations about the problem at hand.
    In the context of (sub)gradient descent, a famous adaptive strategy is the so-called Polyak step-size, which relies on the knowledge of the optimal value  $f_\star$:
    \begin{equation}\label{eq:PolyakStep}
         x_{t+1} = x_t - \gamma_t \nabla f(x_t), \quad \text{with} \ \gamma_t = \frac{f(x_t) - f_\star}{\|\nabla f(x_t)\|^2}.
    \end{equation}
    Polyak steps were originally proposed in~\citet{polyakintroduction} for nonsmooth convex minimization;
    it is also discussed in~\citet{boyd2003subgradient} and a few variants are proposed by, e.g.,~\cite{barre2020complexity,loizou2021stochastic,gower2022cutting} including for stochastic minimization.
    In terms of speed, this strategy (and variants) enjoy similar theoretical convergence properties as those for gradient descent.
    This methods appears to perform quite well in applications where $f_\star$ can be efficiently estimated---see, e.g.,
    ~\cite{hazan2019revisiting} for an adaptation of the method for estimating it online (although~\cite{hazan2019revisiting} contains a number of mistakes,
    the approach can be corrected to achieve the claimed target).
    Therefore, a remaining open question in this context is whether the performances of this method can be improved by incorporating momentum in it.
    A first answer to this question was provided by \citet{barre2020complexity}, although it is not clear that it can match the same convergence properties as optimal first-order methods.

\begin{figure}[ht]
        \centering
        \begin{subfigure}[c]{.9\textwidth}
            \centering
            \includegraphics[trim={0 6.5cm 0 7.5cm},clip, width=\linewidth]{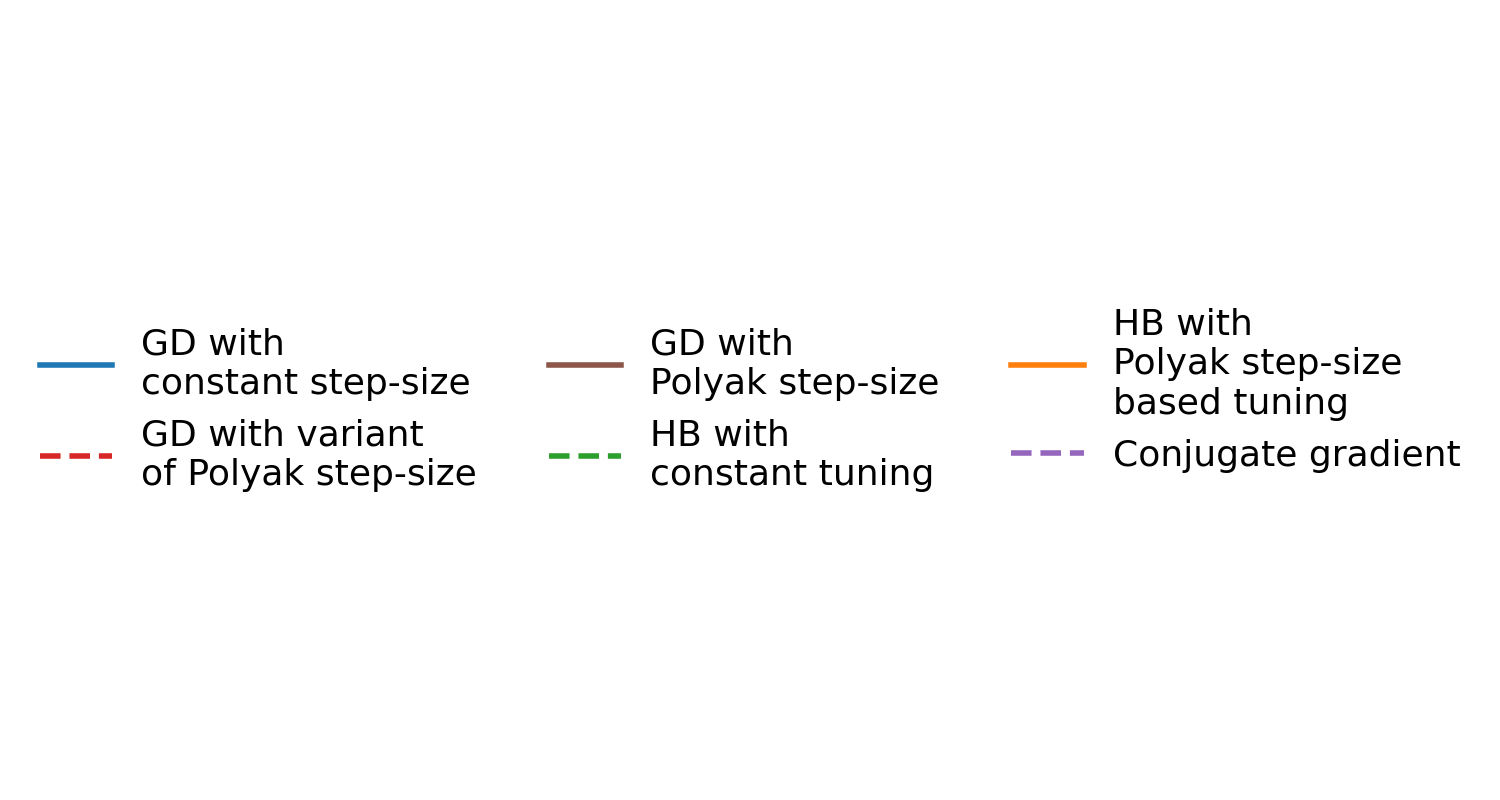}
        \end{subfigure}
        \begin{subfigure}[l]{.75\textwidth}
            \centering
            \includegraphics[width=\linewidth]{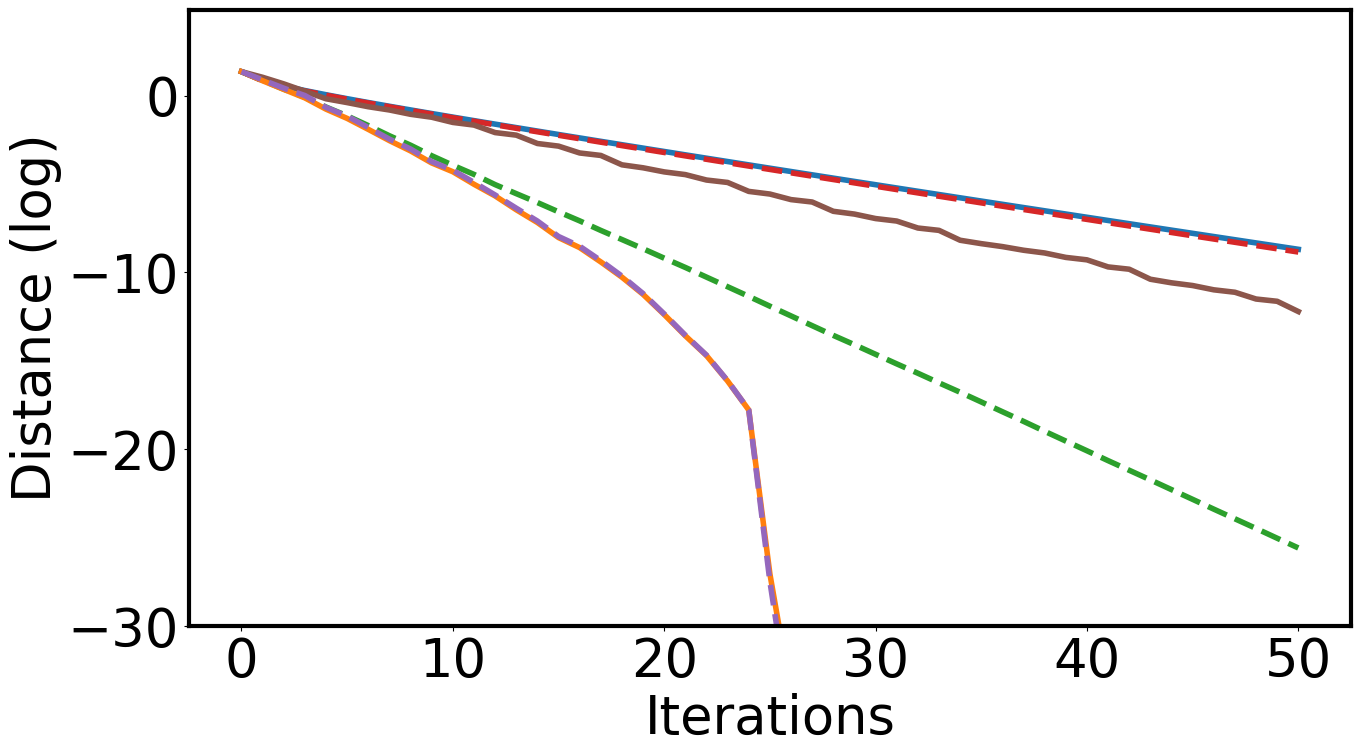}
        \end{subfigure}
        \caption{Comparison in semi-log scale over 50 iterations of different first-order methods applied on a 25-dimensional quadratic objective with condition number $10$. \textbf{GD with constant step-size}, \textbf{GD with Polyak step-size} and \textbf{GD with variant of Polyak step-size} refer to the GD method tuned respectively with the step-size $\gamma = 2 / (L+\mu)$, $\gamma_t = (f(x_t) - f_\star) / \|\nabla f(x_t)\|^2$ and $\gamma_t =   2(f(x_t) - f_\star) / \|\nabla f(x_t)\|^2$. \textbf{HB with constant tuning} is the HB method tuned with constant parameters $\gamma_t = (2 / (\sqrt{L} + \sqrt{\mu}))^2$ and $m_t=((\sqrt{L} - \sqrt{\mu})(\sqrt{L} + \sqrt{\mu}))^2$ while \textbf{HB with Polyak step-size based tuning} refers to Algorithm ~\ref{alg:hb_ps}.}
        \label{fig:quad_short}
    \end{figure}

    \noindent
    In this work, we answer this question for the class of quadratic problems.
    In short, it turns out that the following conjugate gradient-like iterative procedure
    \begin{equation}\label{eq:CG_like}
        \begin{aligned}
            x_{t+1} = \argmin_x \left\{ \| x-\xs\|^2 \text{ s.t. } x \in x_0 + \Span \{ \nabla f(x_0), \nabla f(x_1), \ldots, \nabla f(x_t)\} \right\},
        \end{aligned}
    \end{equation}
    can be rewritten exactly as an Heavy-ball type method whose parameters are chosen adaptively using the value of~$f_\star$.
    This might come as a surprise, as the iteration~\cref{eq:CG_like} might seem impractical due to its formulation relying on the knowledge of $\xs$.
    More precisely,~\cref{eq:CG_like} is exactly equivalent to:
    \begin{equation}\label{eq:PS_momentum}
        \begin{aligned}
            x_{t+1} &= x_t - (1+m_t) \times h_t \nabla f(x_t) + m_t (x_t - x_{t-1}),
        \end{aligned}
    \end{equation}
    with parameters
    \begin{eqnarray}
    \forall t\geq 0, \quad h_t &\triangleq& \frac{2(f(x_t)-f_\star)}{\| \nabla f(x_t) \|^2},\quad   \label{eq:PS_step_tuning}\\
    m_0\triangleq0 \ \  \text{and} \ \  \forall t\geq 1, \quad m_{t} &\triangleq& \frac{-(f(x_{t})-f_\star)\langle \nabla f(x_{t}), \nabla f(x_{t-1}) \rangle}{(f(x_{t-1})-f_\star)\| \nabla f(x_{t}) \|^2 + (f(x_{t})-f_\star)\langle \nabla f(x_{t}), \nabla f(x_{t-1}) \rangle}\ .\label{eq:PS_momentum_tuning}
    \end{eqnarray}
    In~\cref{eq:PS_momentum}, $m_t$ corresponds to the momentum coefficient and $h_t$ to a step-size.
    With the tuning of~\cref{eq:PS_step_tuning}, this step-size is twice the Polyak step-size in~\cref{eq:PolyakStep}.
    This Heavy-ball momentum method with Polyak step-sizes is summarized in~\Cref{alg:hb_ps} and illustrated in~\Cref{fig:quad_short}.
    Due to its equivalence with~\cref{eq:CG_like}, the Heavy-ball method~\cref{eq:PS_momentum} inherits nice advantageous properties of conjugate gradient-type methods, including:
    \begin{enumerate}[label=(\roman*),topsep=0pt,noitemsep]
        \item
    finite-time convergence: the problem~\cref{eq:main-opt} is solved exactly after at most $d$ iterations,
    \item instance optimality: for \textit{all} $H\succcurlyeq 0$, no first-order method satisfying $x_{t+1}\in x_0 + \Span \{ \nabla f(x_0),  \ldots, \nabla f(x_t)\}$ results in a smaller $\|x_t-\xs\|$,
    \item it inherits optimal worst-case convergence rates on quadratic functions.
    \end{enumerate}
    Of course, a few of those points needs to be nuanced in practice due to finite precision arithmetic.
    The equivalence between~\cref{eq:CG_like} and~\cref{eq:PS_momentum} is formally stated in the following theorem.

    \begin{Th}\label{thm:main}
        Let $(x_t)_{t \in \mathbb{N}}$ be the sequence defined by the recursion~\cref{eq:CG_like}, namely such that for any $t$, $x_{t+1}$ is the Euclidean projection of $\xs$ onto the affine subspace $x_0 + \Span \{ \nabla f(x_0), \nabla f(x_1), \ldots, \nabla f(x_t)\}$.
        Then $(x_t)_{t \in \mathbb{N}}$ is the sequence generated by~\Cref{alg:hb_ps}.
    \end{Th}

    \begin{algorithm}[H]
        \caption{Adaptive Heavy-ball algorithm \label{alg:hb_ps}}
        \SetAlgoLined
        \textbf{Input} $T$ and $f: x \mapsto f(x) \triangleq \frac{1}{2} \langle x - \xs,\, H (x - \xs)\rangle + f_\star$\\
        \textbf{Initialize} $x_0\in\mathbb{R}^d$, $m_0=0$ \\
        \For{$t = 0 \cdots T-1$}{
            $h_t = \frac{2(f(x_t)-f_\star)}{\| \nabla f(x_t) \|^2}$\\
            $x_{t+1} = x_t - (1+m_t) \times h_t \nabla f(x_t) + m_t (x_t - x_{t-1}) $\\
            $m_{t+1} = \frac{-(f(x_{t+1})-f_\star)\langle \nabla f(x_{t+1}), \nabla f(x_t) \rangle}{(f(x_t)-f_\star)\| \nabla f(x_{t+1}) \|^2 + (f(x_{t+1})-f_\star)\langle \nabla f(x_{t+1}), \nabla f(x_t) \rangle} $ \;
        }
        \KwResult{$x_T$}
    \end{algorithm}
    
    ~\Cref{thm:main} turns out to be a particular case of a more general result stating that the iterates of any conjugate gradient-type method described with a polynomial $Q$ as
    \begin{equation}\label{eq:CG_like_general}
        x_{t+1} = \mathrm{argmin}_x \left\{ \langle x-\xs,\,  Q(H) (x-\xs)\rangle \text{ s.t. } x \in x_0 + \Span \{ \nabla f(x_0), \ldots, \nabla f(x_t)\} \right\},\tag{Q-minimization}
    \end{equation}
    are equivalently written in terms of an adaptive Heavy-ball iteration.
    In particular,~\cref{eq:CG_like} corresponds to~\cref{eq:CG_like_general} with $Q(x)=1$.
    Similarly, classical conjugate gradient method corresponds to~\cref{eq:CG_like_general} with $Q(x)=x$ (this fact is quite famous,
    see, e.g.,~\cite{nocedal1999numerical}).
    We were surprised not to find this general result written as is in the literature, and we therefore provide it in~\Cref{sec:main_thm}.
    The key point of this work is that the equivalent Heavy-ball reformulation of~\cref{eq:CG_like} can be written in terms of $f_\star$,
    thereby obtaining a momentum-based Polyak step-size.

\paragraph{Notations.} We denote $\preccurlyeq$ the order between symmetric matrices; $\Sp{H}$ the spectrum of the matrix $H$, namely its set of eigenvalues; $\mathbb{R}_d[X]$ the set of polynomials with degree at most $d$.

\subsection{Preliminary material}\label{subsec:preliminary-material}

    \paragraph{Worst-case optimality.} Solving~\cref{eq:main-opt} is a very important problem and several methods have been proposed to achieve this goal.
    They are compared with each other through notions of performance.
    This consists in evaluating the precision of an algorithm over the functions of a given class after a given number $T$ of iterations.
    The main framework is \emph{worst-case analysis} and the precision is the value of a given metric, e.g.\ the distance of the last iterate to optimum $\|x_T-\xs\|$, the function value of the last iterate $f(x_T)-f(\xs)$, or its gradient norm $\|\nabla f(x_T)\|$.
    The \emph{worst-case analysis} framework consists in finding the guarantees of a method that hold for each and every functions of a given class, as for instance the class of $L$-smooth $\mu$-strongly convex quadratic functions described as quadratic functions verifying $\mu I \preceq H \preceq L I$ for given $0 < \mu \leq L$.
    The \textbf{gradient descent (GD)} method characterized by the update
    \begin{equation}
        x_{t+1} = x_t - \gamma_t \nabla f(x_t) \label{eq:gd}
    \end{equation}
    therefore verifies $\|x_{t} - \xs\| = O((\frac{L-\mu}{L+\mu})^{t})$ on all such functions for $\gamma_t = \frac{2}{L+\mu}$.
    Thanks to a relationship with polynomial analysis, \citet{golub1961chebyshev} proved that the \textbf{Chebyshev method}, described as
    \begin{equation}
        x_{t+1} = x_t - \gamma_t \nabla f(x_t) + m_t (x_t - x_{t-1}) \label{eq:hb},
    \end{equation}
    for a well chosen tuning of the parameters $\gamma_t$ and $m_t$ ($m_t = (\frac{\sqrt{L} - \sqrt{\mu}}{\sqrt{L} + \sqrt{\mu}})^2\frac{1 + ((\sqrt{L} - \sqrt{\mu})/(\sqrt{L} + \sqrt{\mu}))^{2(t-1)}}{1 + ((\sqrt{L} - \sqrt{\mu})/(\sqrt{L} + \sqrt{\mu}))^{2(t+1)}}$, $\gamma_t = \frac{2}{L+\mu}(1 + m_t)$), is \emph{worst-case optimal} on this class of function, achieving the guarantee $\|x_{t} - \xs\| = O((\frac{\sqrt{L}-\sqrt{\mu}}{\sqrt{L}+\sqrt{\mu}})^{t})$ (often referred to as ``acceleration'').
    Methods based on this two-term recursion are called ``Heavy-ball'' or ``Polyak momentum''~\citep{polyak1964some}.
    In particular, the stationary regime of the Chebyshev method is the \textbf{Heavy-ball (HB)} method tuned with $m_t = (\frac{\sqrt{L}-\sqrt{\mu}}{\sqrt{L}+\sqrt{\mu}})^{2}$ and $\gamma_t = \frac{2}{L+\mu}(1 + m_t) = (\frac{2}{\sqrt{L} + \sqrt{\mu}})^2$ and achieves the worst-case guarantee
    $\|x_{t} - \xs\| = O( t (\frac{1-\sqrt{\kappa}}{1+\sqrt{\kappa}})^{t})$, close to the optimal one achieved by the Chebyshev method.
    Note that~\cref{eq:hb} is another formulation of~\cref{eq:PS_momentum} where $\gamma_t = (1 + m_t) \times h_t$.
    In all the aforementioned tuning, $h_t$ has the same value: $h_t = \frac{2}{L + \mu}$ (see~\Cref{sec:discussion} for more detailed discussion on this).

    \paragraph{Span of gradients and Krylov subspaces.}
    All methods described above can be defined using a recursion:
    \begin{equation}
        x_{t} = x_0 - \sum_{i=0}^{t-1} \gamma_i^{(t)} \nabla f(x_i) \label{eq:fom}
    \end{equation}
    for some sequence $(\gamma_i^{(t)})_{i\in \llbracket 0, t-1 \rrbracket}$ as suggested in~\citep{nemirovsky1992information, nemirovski1994efficient}.
    Note that the recursion~\cref{eq:fom} can also be explicitly written as $x_{t} = x_0 - H \sum_{i=0}^{t-1} \gamma_i^{(t)} x_i$,
    and therefore, $x_{t} - x_0 \in H \Span(\{x_i\}_{i\in\llbracket 0, t-1 \rrbracket})$.
    We deduce by recursion that $x_t - x_0 \in H \mathcal{K}_t(H, x_0)$ where $\mathcal{K}_t(H, x_0) \triangleq \Span(\{H^i x_0\}_{i\in\llbracket 0, t-1 \rrbracket})$ is called \emph{order-t Krylov subspace generated by $H$ and $x_0$}.
    This creates a link between first-order algorithms and polynomials, summarized in the following lemma (which is implicitely used in~\citet{golub1961chebyshev} and formally stated, e.g., in~\citep[Proposition 4.1]{goujaud2022super}).
    \begin{restatable}{Lemma}{linkalgopoly}\label{lemma:link_algo_poly}.
        Let $f\in\mathcal{C}_{\Lambda}$.
        The iterates $x_t$ satisfy
        \begin{equation}
            x_{t} \in x_0 + \Span\{ \nabla f(x_0),\ldots, \nabla f(x_{t-1}) \} \,, \label{def:first_order_algo}
        \end{equation}
        where $x_0$ is the initial approximation of $\xs$, if and only if there exists a sequence of polynomials $(P_t)_{t\in\mathbb{N}}$, each of degree at most 1 more than the highest degree of all previous polynomials and $P_0$ of degree 0 (hence the degree of $P_t$ is at most $t$), such that
        \begin{equation}
            \vphantom{\sum_i}\forall \;t, \quad x_t - \xs = P_t(H)(x_0-\xs), \quad P_t(0)=1\,. \label{eq:link_polynomial}
        \end{equation}
    \end{restatable}
    Similar to the way we use this technique below, this lemma has already been extensively used to design methods;
    see, e.g.,~\citet[Chapter 1]{d2021acceleration} or the blog post by~\citet{pedregosa2021residual} for gentle introductions to this technique.
    For instance, we can use this technique for optimizing the step-size of the gradient method, or to derive the Chebyshev method, which optimizes the worst-case on the class of smooth and strongly convex quadratic functions~\citep[see][]{fischer2011polynomial}.
    More recently,~\citet{goujaud2022super} used it to derive a method which can take advantage of a possible gap in the spectrum of $H$.
    This approach has also been used for other applications such as accelerated gossip algorithms~\citep{berthier2020accelerated}.

    \paragraph{Adaptive methods.}
    In~\Cref{lemma:link_algo_poly}, while $P_t(H)$ is a polynomial evaluated on the matrix $H$,
    its scalar coefficients might or might not depend on $H$.
    If they depend on $H$, we say that the associated method is adaptive.
    Non-adaptive methods suffer from two main drawbacks:
    (i) they use the same parameters for all the functions within the class of problems, not taking advantage of the observed quantities;
    (ii) the underlying parameters must scale with the function class parameters, and therefore depends on the values of $L$ and $\mu$,
    which are generally difficult to estimate (and actually do not correspond to first-order information, as they rely on the Hessian of the function at hand).
    Ultimately, adaptive methods aim at solving those issues by choosing parameters (step-size, momentum, etc.) on the fly.

    \paragraph{Polyak steps.}

    It is straightforward that a gradient descent update verifies on any convex function $f$ that $\|x_{t+1} - \xs\|^2 \leq \|x_t - \xs\|^2 - 2\gamma_t (f-f_\star) + \gamma_t^2\|\nabla f(x_t)\|^2$.
    ~\cite{polyakintroduction} argues that, based on this inequality, the best guaranteed progression is then achieved for $\gamma_t = \frac{f(x_t) - f_\star}{\|\nabla f(x_t)\|^2}$.
    This choice is called ``Polyak step-size'' and has been studied intensively even recently ~\cite{loizou2021stochastic,gower2022cutting}.

    Other variants of the latter have been proposed.
    For instance~\citep[Variant 1]{barre2020complexity} suggested the step-size $\gamma_t = 2\frac{f(x_t) - f_\star}{\|\nabla f(x_t)\|^2}$.
    This also optimizes the exact progress of one gradient descent update over quadratics realizing a projection of $x_t-\xs$ over the orthogonal subspace of $\nabla f(x_t)$.

    Therefore, the Polyak step-size strategy applied to the gradient descent method achieves the same worst-case guarantee of the well tuned fixed step-size gradient descent method,
    while not relying on Hessian information.
    Moreover, due to its adaptivity to each function, and since generic functions do not look like worst-cases,
    the Polyak step-size strategy applied to the gradient descent method perform very well in practice (See~\Cref{fig:quad_short} and~\citep[Figure 1]{barre2020complexity}),
    sometimes even beating the well tuned non-adaptive Heavy-ball method even if the worst-case guarantees are sorted in a different order.

    \paragraph{Instance-optimality.}
    While optimal worst-case method of the form of~\cref{eq:fom} have been found with predetermined parameters,
    it would be better to find a method under the form of~\cref{eq:fom} that is optimal (for some performance metric),
    not only on worst-case analysis, but on each function individually, taking advantage of the adaptivity of the parameters.
    The well-known conjugate gradient method achieves this goal when the performance metric is the function value of the last iterate.
    The MinRes method attacks the problem minimizing the gradient norm of the last iterate.

    \paragraph{Contributions.}
    In this work we derive iterative methods of the form of~\cref{eq:fom} (which iterates lie in the span of previously observed gradients) that are instance-optimal for a variety of performance metrics.
    All those methods updates are variations of the Heavy-ball two-term recursion~\cref{eq:hb} with only parameters $\gamma_t$ and $m_t$ changing from one method to another.
    Finally, we show (see~\Cref{thm:main}) that for a well chosen yet classical performance metric, this associated method~\Cref{alg:hb_ps} is not relying on second-order information at all (not even $L$ and $\mu$).
    Instead it is using a classical variant of the Polyak step-size $\frac{2(f(x_t) - f_\star)}{\|\nabla f(x_t)\|^2}$, providing an answer to the question ``can we accelerate methods with Polyak step-size?''.

\subsection{Related works}\label{subsec:related-works}

    Polyak step-sizes were proposed in~\citep{polyakintroduction}.
    Despite the dependency on $f_\star$, the Polyak step-size is more studied theoretically and used in practice due to its efficiency when applied to real world problems.
    Recent works~\citep[e.g.][]{loizou2021stochastic, d2021stochastic} argue that this dependency is not a practical issue for many problems which we can assume verify $f_\star=0$ (see~\Cref{sec:experiments}).
    A few variants of the Polyak step-size strategy were proposed by~\citet{barre2020complexity}, including a version incorporating a Nesterov-type momentum~\cite{nesterov1983method}, achieving a worst-case guarantee of $\|x_t - \xs\|^2 = O((1 - 2 (\mu/L)^{2/3})^{2t})$ over the class of (non-necessarily quadratic) $L$-smooth $\mu$-strongly convex functions, thereby improving over previous works on adaptive first-order methods.
    However, the proposed method does not allow to remove the dependency on $L$, and does not achieve the black-box complexity of smooth strongly convex minimization~\cite{nest-book-04}.
    In~\citep{loizou2021stochastic}, the authors study the stochastic Polyak step-size, whereas~\citep{d2021stochastic} applies it to Mirror descent.

    Many alternative adaptive methods have been proposed in the past.
    Among them, let us mention~\citep{barzilai1988two} which introduced the so-called Barzilai-Borwein step-size rule, and the more recent~\citep{malitsky2019adaptive} which developed a step-size policy that adapts to the local geometry with convergence guarantees beyond quadratic minimization.

\section{Main theorem}\label{sec:main_thm}

    This section states and proves~\Cref{thm:main}.
    In short, given a certain function $f$ (characterized by $H$ and $\xs$ here) and a starting point $x_0$,
    we search for an iterative procedure, possibly adaptive, verifying the polynomial-based expression~\cref{eq:link_polynomial} such that $x_t$ converges as fast as possible to $\xs$ for some predefined performance metric.
    Most classical ways to measure the performance of such optimization schemes include the distance to optimum $\|x_t - \xs\|^2$, the function accuracy gap $f(x_t) - f_\star$, the squared gradient norm $\|\nabla f(x_t)\|^2$, and linear combinations of the former.
    Let us abstract those notions by denoting the performance measure of choice by $\langle x_t - \xs,\, Q(H) (x_t - \xs)\rangle$ (with $Q$ a predefined polynomial that is positive on $\mathbb{R}_{>0}$).
    Then, we consider the iterative scheme given by~\eqref{eq:CG_like_general}.
    \begin{equation}\label{eq:q-minimization}
        x_{t+1} = \mathrm{argmin}_x \left\{ \langle x-\xs,\,  Q(H) (x-\xs)\rangle \text{ s.t. } x \in x_0 + \Span \{ \nabla f(x_0), \ldots, \nabla f(x_t)\} \right\} \tag{Q-minimization},
    \end{equation}
    The next theorem provides an explicit instance-optimal method to solve~\cref{eq:q-minimization}.

    \begin{Th}[Main]
        The unique solution to~\cref{eq:q-minimization} is given by the Heavy-ball procedure

        \begin{equation}
            x_{t+1} = x_t - (1 + m_t) \times h_t \nabla f (x_t) + m_t (x_t - x_{t-1})
        \end{equation}

        where

        \begin{equation}
            \left\{
                \begin{array}{ll}
                    h_t & = \frac{\langle x_t - \xs,\, H Q(H) (x_t - \xs)\rangle}{\langle x_t - \xs,\, H^2 Q(H) (x_t - \xs)\rangle}; \\
                    m_t & = \frac{- b_t h_t}{1 + b_t h_t}, \quad \mathrm{~with~} b_t = \frac{\langle x_t - \xs,\, H^2 Q(H) (x_{t-1} - \xs)\rangle}{\langle x_{t-1} - \xs,\, H Q(H) (x_{t-1} - \xs)\rangle}.
                \end{array}
            \right.
            \label{eq:parametrization}
        \end{equation}
    \end{Th}

    Remark that setting $Q(X)=X$ leads to a nice expression of the conjugate gradient method.
    Indeed, setting $Q(X)=X$ corresponds to optimally minimizing the excess loss $f(x_t) - f_\star$.

    As already known, the conjugate gradient method requires the knowledge of H (or an Hessian vector product) to proceed.
    This is also a priori the case for all other choices of $Q(\cdot)$.
    In the case of $Q(X) = 1$, which corresponds to minimizing the distance to the optimum (see~\cref{eq:CG_like}), we can use an alternate writing making use of $f_\star$:
    \begin{equation}
        \left\{
            \begin{array}{ll}
                h_t & = \frac{2(f(x_t) - f_\star)}{\|\nabla f(x_t)\|^2} \\
                m_t & = \frac{- b_t h_t}{1 + b_t h_t}, \quad \mathrm{~with~} b_t = \frac{\left<\nabla f(x_t), \nabla f(x_{t-1}) \right>}{2(f(x_{t-1}) - f_\star)}.
            \end{array}
        \right.
    \end{equation}

\noindent \textit{Proof.}

    \paragraph{Designing methods from the polynomial point of view.}

        As suggested by~\Cref{lemma:link_algo_poly}, we look for an iterative method which can be expressed in the form $x_t - x_\star = P_t(H) (x_0 - x_\star)$, where $P_t$ is a $t^{th}$ degree polynomial with $P_t(0)=1$.
        Furthermore, as we look for an instance-optimal method, the latter polynomial must be instance-specific, and the coefficients of $P_t$ should depend on $H$ (and should describe the iterative procedure~\eqref{eq:CG_like_general}).

        Recalling that $H$ is real symmetric matrix, we denote by $\lambda \in \Sp(H)$ its eigenvalues and by $v_{\lambda}$ th associated orthonormal basis of eigenvectors, leading to $H = \sum_{\lambda\in \Sp(H)} \lambda v_{\lambda} v_{\lambda}^T$.
        The quantity to be minimized can now be written as:

        \begin{align}
            \langle x_t - \xs,\, Q(H) (x_t - \xs)\rangle &~ = \langle x_0 - \xs,\, P_t(H)^T Q(H) P_t(H) (x_0 - \xs)\rangle \\
            &~ = \sum_{\lambda\in \Sp(H)} Q(\lambda) P_t(\lambda)^2 \langle x_0 - \xs, v_{\lambda} \rangle^2 \\
            &~ = \displaystyle\int_{\lambda \in \mathbb{R}^+} P_t(\lambda)^2 ~ \mathrm{d\lambda_Q(\lambda)} \label{eq:distance_x_as_poly}
        \end{align}
        with $\mathrm{\lambda_Q}$ the discrete measure $\sum_{\lambda\in \Sp(H)} Q(\lambda) \langle x_0 - \xs, v_{\lambda} \rangle^2 ~ \delta_{\lambda}$ (we sometimes use the shorthand notation $\int P_t^2 ~ \mathrm{d\lambda_Q}$ for~\cref{eq:distance_x_as_poly} in what follows).
        It is clear that~\cref{eq:distance_x_as_poly} is $0$ if and only if $P_t(\lambda) = 0$ for all $\lambda\in \Sp(H)$.
        As a consequence, we conclude that (i) choosing the right sequence of polynomials leads to convergence in exactly $|\Sp(H)|$ iterations, and (ii) $\langle P^{(1)}, P^{(2)} \rangle_{Q} \triangleq \int P^{(1)}P^{(2)} ~ \mathrm{d\lambda_Q}$ is an inner product on $\mathbb{R}_{|\Sp(H)|-1}[X]$.
        We therefore want to solve
        \begin{equation}
            \left\{
            \begin{array}{ll}
                \underset{P_t \in \mathbb{R}_t[X]}{\mathrm{minimize~}} & \| P_t \|_Q^2 \\
                \mathrm{subject~to~} & P_t(0)=1
            \end{array}
            \right.
            \label{eq:main_problem}
        \end{equation}
        for any $t \leq |\Sp(H)| - 1$ where $\| P \|_Q^2 \triangleq \langle P, P \rangle_{Q} = \int P^2 ~ \mathrm{d\lambda_Q}$ denotes the underlying norm of the inner product $\langle \cdot, \cdot \rangle_{Q}$.
        For $t \geq |\Sp(H)|$, we consider instead $P_t$ as a multiple of the polynomial $\prod_{\lambda\in \Sp(H)} (X - \lambda)$ in~$X$.
        The next steps are somewhat standard and follow a classical pattern for solving~\cref{eq:main_problem} (see, e.g.\ \citet{berthier2020accelerated} and the references therein).

    \paragraph{From minimal norm to orthogonality.}

        The solution to~\cref{eq:main_problem} is the projection of the polynomial $0$ over the affine space
        $\left\{ P\in\mathbb{R}_t[X] \mid P(0)=1 \right\}$ with respect to the inner product $\langle \cdot, \cdot\rangle_Q$.
        A necessary and sufficient condition for $P$ to be the solution of problem~\cref{eq:main_problem} is therefore
        to verify $\langle 0-P, \Delta P \rangle_Q=0$ for any $\Delta P$ in the vectorial subspace
        $\left\{ P\in\mathbb{R}_t[X] \mid P(0)=0 \right\}=X\mathbb{R}_{t-1}[X]$.
        Hence $P_t$ solves problem~\cref{eq:main_problem} iff
        \begin{equation}
            \langle P_t, XR \rangle_Q=0, \forall R\in\mathbb{R}_{t-1}[X].
        \end{equation}
        Note however, that for any $(P, R)\in\mathbb{R}[X]^2$,
        \begin{align*}
            \langle P, XR \rangle_Q
            & = \displaystyle\int_{\lambda\in \mathbb{R}^+} P(\lambda) \times \lambda R(\lambda) ~ \mathrm{d\lambda_Q(\lambda)} \\
            & = \displaystyle\int_{\lambda\in \mathbb{R}^+} P(\lambda) \times R(\lambda) ~ \mathrm{d\lambda_{XQ}(\lambda)} \\
            & \triangleq \langle P, R \rangle_{XQ}
        \end{align*}
        with $\mathrm{d\lambda_{XQ}(\lambda)} \triangleq \mathrm{\lambda d\lambda_{Q}(\lambda)}
        = \displaystyle \sum_{\lambda\in \Sp(H)} \lambda Q(\lambda) \times \langle x_0 - \xs, v_{\lambda} \rangle^2 ~ \delta_{\lambda}$.
        Using the latter inner product, the condition for $P_t$ to be the solution to problem~\cref{eq:main_problem} becomes:
        \begin{equation}
            P_t \in \mathbb{R}_{t-1}[X]^{\perp_{XQ}}
        \end{equation}
        Hence, $(P_t)_{t\in\mathbb{N}}$ is a family of orthogonal polynomials for the inner product $\langle \cdot, \cdot \rangle_{XQ}$.

    \paragraph{From orthogonality to recursion.}
        We now focus on finding an explicit expression for the polynomials $P_t$.
        As for all families of orthogonal polynomials, $(P_t)_{t\in\mathbb{N}}$ can be obtained through a
        two-term recursion of the form:
        \begin{equation}
            P_{t+1}(X) = (a_t X + b_t) P_t(X) + c_t P_{t-1}(X), ~\text{ for some } (a_t, b_t, c_t)\in\mathbb{R}^3, \label{eq:3_terms_recursion}
        \end{equation}
        which is easy to verify by induction.
        Our goal is to find $a_t$, $b_t$ and $c_t$.
        First, notice that $(a_t X + b_t) P_t(X) + c_t P_{t-1}(X)$ is orthogonal to $\mathbb{R}_{t-2}[X]$ independently of the values of $a_t$, $b_t$ and $c_t$.
        Those three coefficients can be found via the following  three conditions: (i)
        $\langle P_{t+1}, P_t \rangle_{XQ} = 0$, (ii)
        $\langle P_{t+1}, P_{t-1} \rangle_{XQ} = 0$, and (iii) $P_{t+1}(0)=1$.

        More precisely, it is clear that $a_t\neq 0$ for $P_{t+1}$ to be of degree $t+1$.
        Therefore, one can factorize by $a_t$.
        Reparametrizing~\cref{eq:3_terms_recursion}, one can write
        \begin{equation*}
            P_{t+1}(X) = \frac{(\tilde{a}_t - X) P_t(X) + \tilde{b}_t P_{t-1}(X)}{\tilde{c}_t}, ~\text{with } (\tilde{a}_t, \tilde{b}_t, \tilde{c}_t)\in\mathbb{R}^3.
        \end{equation*}

        Moreover, evaluation at $X=0$ gives $\frac{\tilde{a}_t + \tilde{b}_t}{\tilde{c}_t}=1$, thereby enforcing $\tilde{c}_t = \tilde{a}_t + \tilde{b}_t$.
        It remains to verify the two orthogonality conditions (independent of $\tilde{c}_t$):

        \begin{equation*}
            \begin{array}{lll}
                \tilde{a}_t \langle P_t, P_t \rangle_{XQ} & + \tilde{b}_t \langle P_{t-1}, P_t \rangle_{XQ} & = \langle X P_t(X), P_t(X) \rangle_{XQ}, \\
                \tilde{a}_t \langle P_t, P_{t-1} \rangle_{XQ} & + \tilde{b}_t \langle P_{t-1}, P_{t-1} \rangle_{XQ} & = \langle X P_t(X), P_{t-1}(X) \rangle_{XQ}.
            \end{array}
        \end{equation*}
        Note that this system of equations is decoupled since $\langle P_{t-1}, P_t \rangle_{XQ}=0$, and we finally arrive to
        \begin{equation}
            P_{t+1}(X) = \frac{(\tilde{a}_t - X) P_t(X) + \tilde{b}_t P_{t-1}(X)}{\tilde{a}_t + \tilde{b}_t}, \label{eq:final_poly_recursion}
        \end{equation}
        with
        \begin{equation}
            \left\{
                \begin{array}{ll}
                    \tilde{a}_t & = \frac{\langle X P_t(X), P_t(X) \rangle_{XQ}}{\langle P_t, P_t \rangle_{XQ}}, \\
                    \tilde{b}_t & = \frac{\langle X P_t(X), P_{t-1}(X) \rangle_{XQ}}{\langle P_{t-1}, P_{t-1} \rangle_{XQ}}.
                \end{array}
            \right.
            \label{eq:an_and_bn_values}
        \end{equation}

    \paragraph{From a polynomial recursion to an iterative optimization method.}

        For reaching the final desired result, we simply multiply~\cref{eq:final_poly_recursion} (evaluated in $H$) by $x_0 - \xs$:
        \begin{align}
            x_{t+1} - \xs &~ = \frac{\tilde{a}_t (x_t - \xs) - H (x_t - \xs) + \tilde{b}_t (x_{t-1} - \xs)}{\tilde{a}_t + \tilde{b}_t}, \nonumber \\
            &~ = x_t - \xs - \frac{1}{\tilde{a}_t + \tilde{b}_t} \nabla f (x_t) + \frac{- \tilde{b}_t}{\tilde{a}_t + \tilde{b}_t} (x_t - x_{t-1}), \nonumber
        \end{align}
    thereby reaching the desired
    \begin{equation}
        x_{t+1} ~ = x_t - h_t \nabla f (x_t) + m_t (x_t - x_{t-1})
            \label{eq:HB_formulation}
    \end{equation}
        with
        \begin{equation}
            h_t = \frac{1}{\tilde{a}_t + \tilde{b}_t}, \quad \text{and} \quad m_t = \frac{- \tilde{b}_t}{\tilde{a}_t + \tilde{b}_t}.
            \label{eq:HB_tuning}
        \end{equation}
    From~\cref{eq:HB_formulation}, we recognize an Heavy-ball method with some variable step-size $h_t$ and momentum term $m_t$, thereby concluding the proof. $\hfill\blacksquare$

\begin{Rem}[Step-size parametrization.]
    While the $\gamma_t$ plays a different role in~\cref{eq:hb} and~\cref{eq:gd}, they both usually are called ``step-size'' by default.
    But we noticed that both in Chebyshev method and the Heavy-ball method (optimally tuned), $h_t = \frac{\gamma_t}{1 + m_t}$ is exactly $\frac{2}{L + \mu}$, value of the optimal step-size for gradient descent.
    In~\cref{eq:PS_step_tuning}, we notice again that the value of $h_t$ is the optimal step-size for a single step of gradient descent.
    For this reason, we believe that the natural parametrization of the Heavy-ball methods should be $x_{t+1} = x_t - (1 + m_t) \times h_t \nabla f(x_t) + m_t (x_t - x_{t-1})$ and that $h_t$ should be referred to as the ``natural'' step-size.
    Indeed, when one thinks of the Heavy-ball method with Polyak step-sizes, they would set $\gamma_t$ to the Polyak step-size, not $h_t = \frac{\gamma_t}{1 + m_t}$.
    We therefore provide a novel view on what should be tested.
\end{Rem}

\section{Numerical experiments}\label{sec:experiments}

\begin{figure}[H]
        \centering
        \begin{subfigure}[c]{.9\textwidth}
            \centering
            \includegraphics[trim={0 6.5cm 0 7.5cm},clip, width=\linewidth]{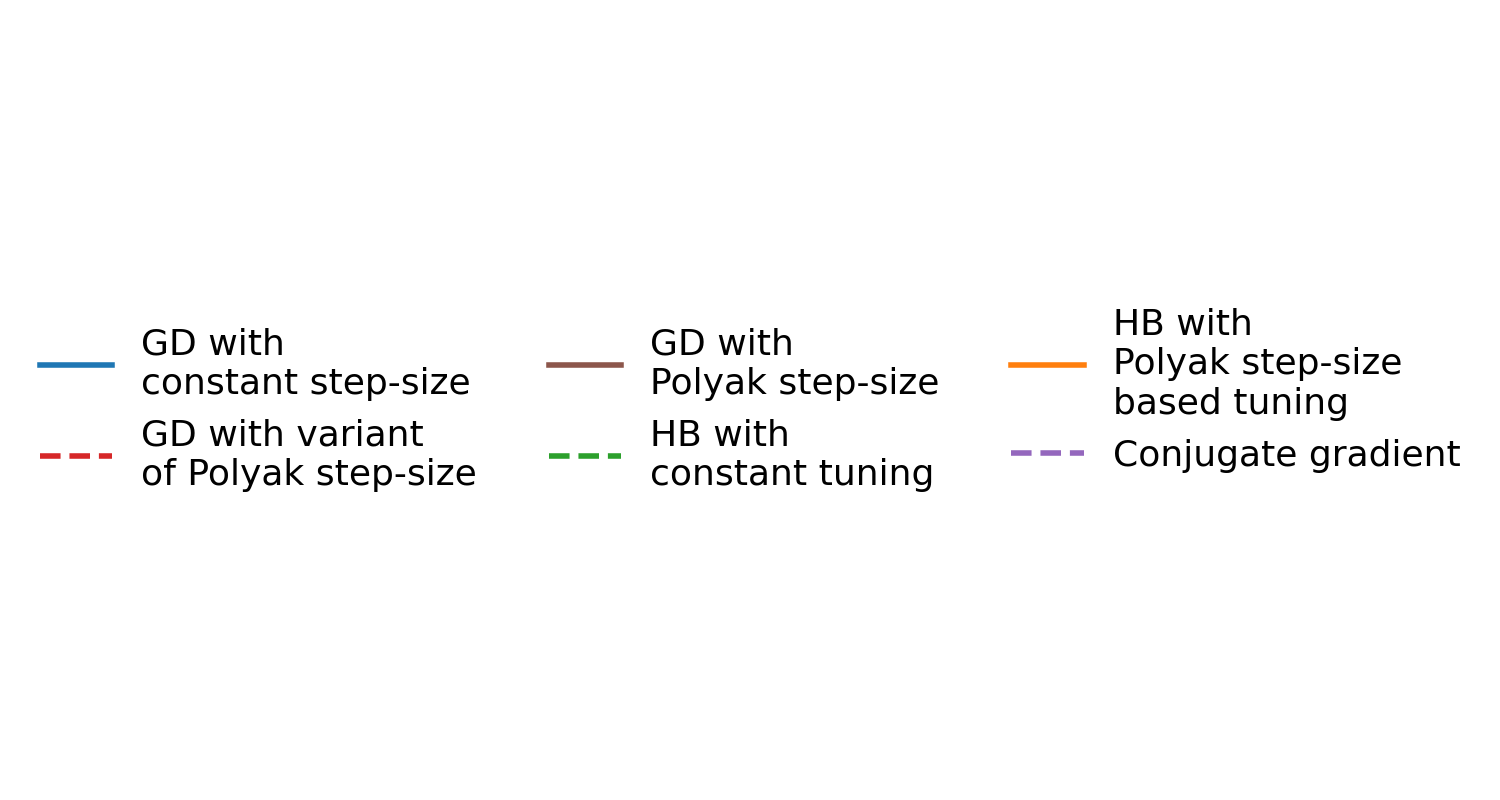}
        \end{subfigure}
        \begin{subfigure}[l]{.45\textwidth}
            \centering
            \includegraphics[width=\linewidth]{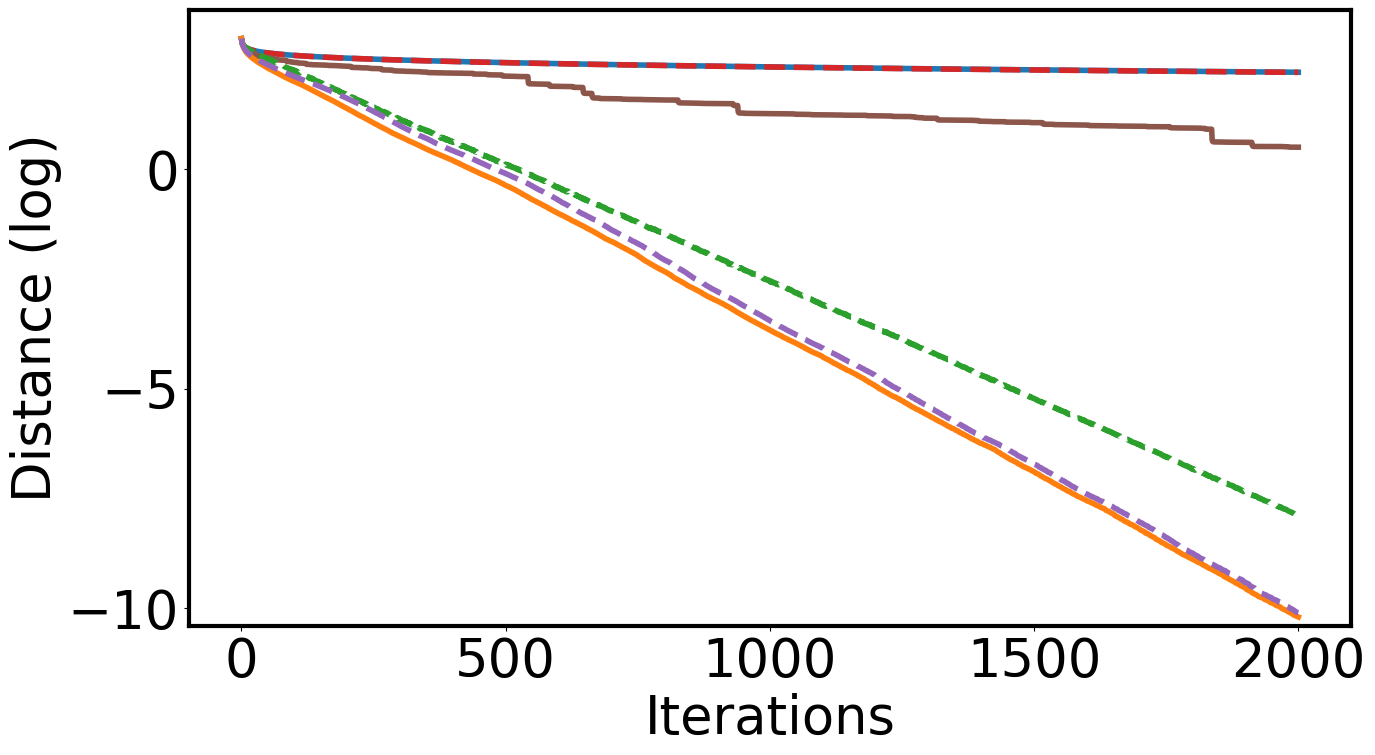}
            \caption{Comparison of distances to optimum}
            \label{fig:quad_long_distances}
        \end{subfigure}
        \hspace{.2cm}
        \begin{subfigure}[r]{.45\textwidth}
            \centering
            \includegraphics[width=\linewidth]{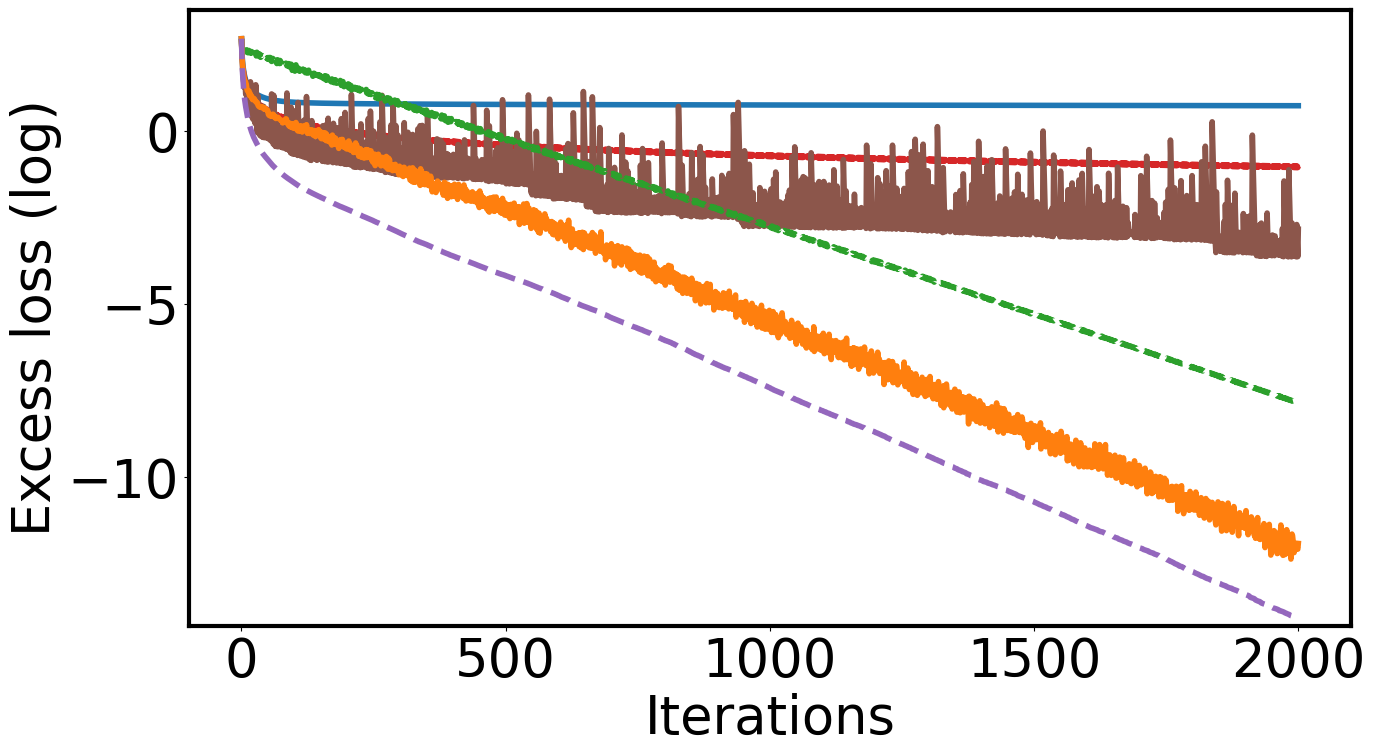}
            \caption{Comparison of excess losses}
            \label{fig:quad_long_losses}
        \end{subfigure}
        \caption{Comparison in semi-log scale over 2000 iterations of different first-order methods applied on a 1000-dimensional quadratic objective with condition number $10^5$. \textbf{GD with constant step-size}, \textbf{GD with Polyak step-size} and \textbf{GD with variant of Polyak step-size} refer to the GD method tuned respectively with the step-size $\gamma = 2 / (L+\mu)$, $\gamma_t = (f(x_t) - f_\star) / \|\nabla f(x_t)\|^2$ and $\gamma_t =   2(f(x_t) - f_\star) / \|\nabla f(x_t)\|^2$. \textbf{HB with constant tuning} is the HB method tuned with constant parameters $\gamma_t = (2 / (\sqrt{L} + \sqrt{\mu}))^2$ and $m_t=((\sqrt{L} - \sqrt{\mu})(\sqrt{L} + \sqrt{\mu}))^2$ while \textbf{HB with Polyak step-size based tuning} refers to Algorithm ~\ref{alg:hb_ps}.}
        \label{fig:quad_long}
    \end{figure}

    In this section, we compare gradient descent, Heavy-ball and conjugate gradient method in an adaptive setting or not. \Cref{fig:quad_long} shows the performance of all these methods on a quadratic objective with known minimal value $f_\star$.
    The hessian of this quadratic objective has been generating from a sequence of eigenvalues with geometric increase, and a random orthogonal transformation.
    The difference between~\Cref{fig:quad_short} and~\Cref{fig:quad_long} is the dimension of the problem as well as the condition number of the objective function.
    Due to finite precision arithmetic, the finite-time convergence is not visible when the condition number is too large.
    However both figures show that our method and the conjugate gradient algorithm behave similarly and faster than the other methods. The code can be found on this~\href{https://github.com/bgoujaud/Heavy-ball_polyak_steps}{GitHub repository}.

\section{Concluding remarks and discussion}\label{sec:discussion}
    
    Polyak step-sizes are known for their general good working performances when the optimal value to the optimization problem at hand is known. The question of whether Polyak step-sizes can be used together with momentum for obtaining accelerated first-order methods appears to be an open question~\citep{barre2020complexity}, which we answer in the simpler case of convex quadratic minimization. In this context, we argue that not only this tuning works well, but also it pops up naturally when investigating instance-optimal first-order iterative methods. Furthermore, we believe it is a necessary step for being able to understand more general optimization settings beyond quadratics. As our method does not seem to work well beyond quadratics, we leave further investigations on this topic for future work.
    
    Among our competitors, we note that the celebrated conjugate gradient (CG) method is another instance-optimal algorithm for quadratics. Whereas our method minimizes the distance to the solution at each iteration, CG is instance-optimal for minimizing function values at each iteration. Perhaps interestingly, the two methods appeared to behave similarly in our numerical experiments. That being said, the main practical differences between the two methods are that CG Heavy-ball-like formulation naturally relies on higher order information while Polyak step-sizes do require knowledge of $f_\star$. In typical optimization problems, this value is not known. However, there are a few settings where this value is actually well-known, typically when $f_\star=0$ generically (in machine learning, this setting is known as the ``interpolation'' regime; an alternative could be to use Polyak-steps as a competitor to MinRes). Finally, let us mention that although a few generalization of CG, often referred to as nonlinear conjugate gradient, were studied in the literature (see, e.g.,~\citep{bonnans2006numerical,nocedal1999numerical,hager2006survey}).

\section*{Acknowledgement}
    We would to thank Raphael Berthier for his insightful feedbacks.
    The work of B. Goujaud and A. Dieuleveut is partially supported by ANR-19-CHIA-0002-01/chaire SCAI, and Hi!Paris.
    A. Taylor acknowledges support from the European Research Council (grant SEQUOIA 724063).
    This work was partly funded by the French government under management
    of Agence Nationale de la Recherche as part of the ``Investissements d’avenir'' program,
    reference ANR-19-P3IA-0001 (PRAIRIE 3IA Institute).

\bibliographystyle{plainnat}
\bibliography{references}

\end{document}